\documentclass{amsart}
\usepackage{amsmath,amsthm}
\usepackage{amsfonts,amssymb}

\usepackage{enumerate}

\newtheorem{thm}{Theorem}[section]
\newtheorem{cor}[thm]{Corollary}
\newtheorem{lem}[thm]{Lemma}
\newtheorem{prop}[thm]{Proposition}

\newtheorem{defn}[thm]{Definition}

\theoremstyle{remark}
\newtheorem{rem}{Remark}[section]

 \def\tr{{\triangle}}
 \def\ll{\lesssim}

\def\sph{\mathbb{S}^{d-1}}
\def\f{\frac}

 \def\a{{\alpha}}

 \def\t{{\theta}}
 \def\l{{\lambda}}
 
 \def\o{{\omega}}
 \def\s{{\sigma}}
 \def\la{{\langle}}
 \def\ra{{\rangle}}

 \def\CD{{\mathcal D}}
 
 \def\CH{{\mathcal H}}

 \def\CV{{\mathcal V}}
 \def\CW{{\mathcal W}}
 \def\BB{{\mathbb B}}

 \def\NN{{\mathbb N}}

 \def\RR{{\mathbb R}}
  \def\SS{{\mathbb S}}

        \def\vi{\varphi}
        \def\p{\partial}

\newcommand{\wt}{\widetilde}
\newcommand{\wh}{\widehat}

\begin{document}

\title[Polynomial approximation in Sobolev spaces]
{Polynomial Approximation in Sobolev Spaces on the Unit Sphere and the Unit Ball}
\author{Feng Dai}
\address{Department of Mathematical and Statistical Sciences\\
University of Alberta\\, Edmonton, Alberta T6G 2G1, Canada.}
\email{dfeng@math.ualberta.ca}
\author{Yuan Xu}
\address{Department of Mathematics\\ University of Oregon\\
    Eugene, Oregon 97403-1222.}\email{yuan@math.uoregon.edu}
\date{\today}
\thanks{Work of the  first  author  was  supported  in part
 by the NSERC grant of Canada.}

\begin{abstract}
This work is a continuation of the recent study by the authors on approximation theory
over the sphere and the ball. The main results define new Sobolev spaces on these
domains and study polynomial approximations for functions in these spaces, including
simultaneous approximation by polynomials and relation between best approximation
to a function and to its derivatives.
\end{abstract}

\maketitle

\section{Introduction}
\setcounter{equation}{0}

In a recent work \cite{DaX}, the authors defined new moduli of smoothness
and K-functionals on the unit sphere $\SS^{d-1} = \{x\in \RR^d: \|x\| =1\}$ and
the unit ball $\BB^d=\{x\in\RR^d: \|x\| \le 1\} $ of $\RR^d$, where $\|x\|$
denotes the usual Euclidean norm, and used them to characterize the
best approximation by polynomials. This work is a continuation of \cite{DaX}
and studies polynomial approximation in Sobolev spaces.

The new modulus of smoothness on the sphere is defined in terms of forward
differences in the angle $\t_{i,j}$ of the polar coordinates on the $(x_i,x_j)$
planes, and it is essentially the maximum over all possible angles of the moduli
of smoothness of one variable in $\t_{i,j}$. There are $\binom{d}{2}$ such angles,
which are clearly redundant as a coordinate system. Nevertheless, our new
definition effectively reduces a large part of problems in approximation theory
on $\SS^d$ to problems on $\SS^1$, which allows us to tap into the rich
resources of trigonometric approximation theory for ideas and tools and adopt
them for problems on the sphere. These same angles also become indispensable
for our new definition of moduli of smoothness on the unit ball $\BB^d$. In fact,
our moduli on $\BB^d$ are defined as the maximum of moduli of smoothness of
one variable in these angles and of one additional term that takes care of the
boundary behavior. We had two ways to define the additional term, the
first one is deduced from the results on the sphere and the second one is the
direct extension of the Ditzian-Totik modulus of smoothness on $[-1,1]$, both
of which capture the boundary behavior of the unit ball and permit both direct and
inverse theorem for the best approximation. For $d=1$, approximation by polynomials
on $\BB^1 = [-1,1]$ is often deduced from approximation by trigonometric
polynomials on the circle $\SS^1$ by projecting even functions and their approximations
on $\SS^1$ to $[-1,1]$. This procedure can be adopted to higher dimension
by projection functions on $\SS^d$ onto $\BB^d$, and this is how our first
modulus of smoothness on $\BB^d$ was defined.  It should be mentioned that
our new moduli of smoothness on the sphere and on the ball are computable;
in fact, the computation is not much harder than what is needed for computing
classical modulus of smoothness of one variable. A number of examples were
given in \cite{DaX}.

In the present paper we continue the work in this direction, study best approximation 
of functions and their derivatives. On the sphere, our result will be given in terms 
of differential operators $D_{i,j} = x_i \partial_j - x_j \partial_i$, which can be 
identified as partial derivatives with respect to $\t_{i,j}$. We shall define Sobolev 
spaces and Lipschitz  spaces in terms of $D_{i,j}$ on these two domains and study 
approximation by polynomials in these spaces, including simultaneous approximation
of functions and their derivatives. The study is motivated by a question of 
Kendall Atkinson (cf. \cite{ACH}) in numerical solution of Poisson equation. 

The paper is organized as follows. The main 
results in \cite{DaX} on the unit sphere will be recalled in Section 2 and the new 
results on the sphere will be developed in Section 3. The results in \cite{DaX} on 
the unit ball will be recollected and further clarified in Section 4. Finally, the new 
results on the ball are developed in Section 5.

Throughout this paper we denote by $c, c_1, c_2, ...$ generic constants that
may depend on fixed parameters and their values may vary from line to line.
We write $A \ll B$ if $A \le c B$ and $A \sim B$ if $A \ll B$ and $B \ll A$.

\section{Polynomial approximation on the sphere: Recent Progress}
\setcounter{equation}{0}

In this section we recall recent progress on polynomial approximation on
the sphere as developed in \cite{DaX}.

Let $L^p(\sph)$ be the $L^p$-space with respect to  the usual Lebesgue measure
$d\s$ on $\sph$ with norm denoted by $\| \cdot \|_{p} := \|\cdot\|_{L^p(\sph)}$ for
$1 \le p \le \infty$, where for $p = \infty$, we replace $L^\infty$ by $C(\sph)$,
the space of continuous functions on $\sph$ with the uniform norm.

\subsection{Polynomial spaces and spherical harmonics}
We denote by $\Pi_n^d$  the space of polynomials of total degree $n$ in $d$ variables,
and by $\Pi_n(\sph):=\Pi_n^d \vert_{\sph}$ the space of all  polynomials in $\Pi_n^d$
restricted on $\sph$. In the following we shall write $\Pi_n^d$ for $\Pi_n^d(\sph)$
whenever it causes no confusing. The quantity of best approximation is then defined by
\begin{equation}\label{bestapp}
    E_n (f)_{p} : = \inf_{ g  \in \Pi_{n-1}^d} \| f - g \|_{p}, \qquad 1 \le p \le \infty.
\end{equation}

Let $\CH_n^d$ denote the space of spherical harmonics of degree
$n$ on $\sph$, which are the restrictions of homogeneous harmonic
polynomials to $\sph$.  Let $\Delta_0$ be the Laplace-Beltrami
operator on the sphere, defined by
\begin{equation} \label{Laplace-Betrami}
    \Delta_0 f(x) : = \Delta \left[ f \left(\frac{y}{\|y\|}\right) \right](x),\   \   \    \ x \in\SS^{d-1},
\end{equation}
where $\Delta: = \frac{\partial^2}{\partial x_1^2}+ \ldots +\frac{\partial^2}{\partial x_d^2}$
is the usual Laplace operator and it acts on the variables $y$. Then the
spherical harmonics are the eigenfunctions of $\Delta_0$,
\begin{equation} \label{Laplace-Betrami-eigen}
    \Delta_0 Y = - n (n+d-2) Y, \qquad Y \in \CH_n^d.
\end{equation}
The reproducing kernel of the space $\CH_n^d$ in $L^2(\sph)$ is given by
the zonal harmonic
\begin{equation}\label{zonal}
    Z_{n,d}(x,y) := \f{n + \l}{\l} C_n^\l(\la x,y\ra), \qquad \l = \f{d-2}2,
\end{equation}
where $\la x, y\ra$ denotes the Euclidean dot product of $x, y\in\RR^d$. Let $C_n^\l$
be the Gegenbauer polynomial with index $\l$, normalized by $C_n^\l(1) =
\binom{n+2\l-1}{n}$. Let $\eta$ be a $C^\infty$-function on $[0,\infty)$ with the
properties that $\eta(x)=1$ for $0\leq x\leq 1$ and $\eta(x)=0$ for $x\ge 2$. We define
\begin{equation}\label{Vnf}
V_nf(x): =\int_{\sph} f(y) K_n(\la x, y\ra)\, d\s(y), \quad
x\in\sph, \quad n=1,2,\cdots,
\end{equation}
where
$$
K_n(t) :=\sum_{k=0}^{2n} \eta\left(\f kn\right)
\f{k+\l}{\l}C_k^{\l}(t),
    \quad  t\in[-1,1].
$$
Then $V_n f \in \Pi_{2n}^d$, $V_ n f = f $ for all $f \in \Pi_n^d$, and for $f \in L^p(\sph)$
or $C(\sph)$,
\begin{equation} \label{Vnf-best}
     \| V_n f - f\|_p \le c E_n(f)_p,  \quad 1 \le p \le \infty.
\end{equation}

\subsection{A class of differential operators on $\SS^{d-1}$}
One of the  main tools in our study is  a class of differential
operators $D_{i,j}$, $1\leq i\neq j\leq d$, that commute with the
Laplace-Beltrami operator. Let $SO(d)$ denote the group of
rotations on $\RR^d$, and let $e_1, \cdots, e_d$ denote the
standard orthogonal basis in $\RR^d$. For $1\leq i \neq j \leq d$
and $t \in \RR$, we denote by $Q_{i,j,t}$ the  rotation by the
angle $t$ in the $(x_i, x_j)$-plane, oriented such that the
rotation from the vector $e_i$ to the vector $e_j$ is assumed to
be positive. For example, the action of the rotation $Q_{1,2,
t}\in SO(d)$ is given  by
\begin{align}\label{EulerAngle1}
Q_{1,2,t}(x_1,\ldots, x_d) = & (x_1 \cos t -x_2 \sin t, x_1 \sin t
+ x_2 \cos t, x_3,
   \cdots, x_d) \\
  = & (s \cos (\phi + t), s \sin (\phi +t), x_3, \cdots, x_d), \notag
\end{align}
where $(x_1,x_2) = s (\cos \phi, \sin \phi)$, and other $Q_{i,j,t}$ are defined likewise.

To each $Q \in SO(d)$ corresponds an operator $L(Q)$ in the space $L^2(\SS^{d-1})$,
defined by $L(Q) f(x):= f(Q^{-1} x)$ for $x \in \SS^{d-1}$, which is a group representation
of $SO(d)$. The infinitesimal operator of $L(Q_{i,j,t})$ has the form
\begin{equation}\label{partial_ij}
 D_{i,j}: = \frac{\partial}{\partial t} \left[ L(Q_{i,j,t})\right] \Big \vert_{t =0}
       = x_j \frac{\partial}{\partial x_i}   - x_i \frac{\partial}{\partial x_j}, \qquad
        1 \le i < j \le d.
\end{equation}
In particular, it is easy to verify that, taking $(i,j) = (1,2)$ as an example,
\begin{equation}\label{partial_ij2}
 D_{1,2}^r f(x) = \left(-\frac{\partial}{\partial \phi}\right)^r
  f(s \cos \phi, s \sin \phi, x_3,\ldots,x_d).
 \end{equation}
The following useful observation, which asserts that  $D_{i,j}^r f$ is independent of smooth
extensions of $f$, is a simple consequence of \eqref{partial_ij2}:

\begin{prop}
Let $x_0\in \sph$, and let  $F$ and $G$ be  two smooth functions on an open
neighborhood $U\subset \RR^d$ of $x_0$ which coincide on $U\cap \sph$. Then
$D_{i,j}^r F(x_0)=D_{i,j}^r G (x_0).$
\end{prop}

The operators $D_{i,j}$ are  connected to the usual tangential
partial derivatives according to the following formula
\cite[(3.15)]{DaX}:
\begin{equation} \label{partial-Dij}
   \frac{\partial}{\partial x_j} \left[ f\left(\frac{x}{\|x\|}\right)\right]_{\|x\|=1}
      = \partial_j f - x_j \sum_{i=1}^d x_i \partial_i f =
          - \sum_{\{i: 1\leq i\neq j\leq d\}} x_i D_{i,j} f.
\end{equation}
The operators $D_{i,j}$ are also closely related to the Laplace-Beltrami operator $\Delta_0$.
In fact, $\Delta_0$ satisfies the following decomposition \cite[(2.6)]{DaX},
\begin{equation} \label{Laplace-Betrami2}
  \Delta_0 = \sum_{1 \le i < j \le d} D_{i,j}^2.
\end{equation}
Furthermore, each operator $D_{i,j}$ in this decomposition commutes with $\Delta_0$.
In particular, by \eqref{Laplace-Betrami-eigen}, this implies that the  spaces of spherical
harmonics on $\sph$ are  invariant under $D_{i,j}$.

\subsection{Moduli of smoothness and K-functionals on $\SS^{d-1}$}
For  $1 \le i <j \le d$ and $\t \in [-\pi, \pi]$, we define the
$r$th difference operator $ \tr_{i,j,\t}^r$ by
$$\tr_{i,j,\t}^r :=(I-T_{Q_{i,j,\t}})^r = \sum_{k=0}^r (-1)^k \binom { r}
k  T_{Q_{i,j,k\t}},$$ where
 $T_Q f(x): =f(Q x)$ for $Q\in SO(d)$.  This differential operator can be expressed in terms
of the usual forward difference as, for example, for $(i,j) = 1(,2)$,
\begin{equation} \label{Delta_ij}
\tr_{1,2, \t}^r f(x) =  \overrightarrow{\tr}_\t^r
   f \left (x_1 \cos(\cdot)- x_2 \sin (\cdot), x_1 \sin(\cdot) + x_2 \cos (\cdot),
        x_3,\ldots,x_d \right),
\end{equation}
where $\overrightarrow{\tr}_\t^r$ is acted on the variable
$(\cdot)$, and is evaluated at $t=0$. The following new modulus of
smoothness  was  recently introduced  in \cite[Definition 2.2]{DaX}:

\begin{defn}
For $r \in \NN$, $t > 0$, and $f \in L^p(\SS^{d-1})$, $1 \le p <
\infty$, or $f \in C(\SS^{d-1})$ for $p = \infty$, define
\begin{equation} \label{eq:modulus}
 \o_r (f,t)_p : = \sup_{ |\t| \le t} \max_{1 \le i < j \le d}
          \left \|\Delta_{{i,j,\t}}^r f \right \|_p.
\end{equation}
\end{defn}

This modulus of smoothness enjoys most of the properties of classical moduli
of smoothness (\cite[Proposition 2.7]{DaX}), and it permits both direct and
inverse theorems (\cite [Theorem 3.4]{DaX}):

\begin{thm} \label{thm:best}
For $f \in L^p(\sph)$ if $1 \le p < \infty$ and $f \in C(\sph)$ if $p
=\infty$,
\begin{equation} \label{Jackson}
       E_n(f)_p \le c \, \o_r(f,n^{-1})_p, \qquad 1 \le p \le \infty.
\end{equation}
On the other hand,
\begin{equation} \label{inverse}
       \o_r(f,n^{-1})_p \le c\, n^{-r} \sum_{k=1}^n k^{r-1} E_{k-1}(f)_p, \qquad 1 \le p \le \infty.
\end{equation}
\end{thm}

A new K-functional on $\sph$ was defined in terms of the differential operators $D_{i,j}$
in \cite[Definition 2.4]{DaX}:

\begin{defn}   \label{def:K-func}
For $r \in \NN_0$ and $t \ge 0$,
\begin{equation} \label{eq:K-func-sphere}
   K_r(f,t)_p : = \inf_{g \in C^r(\sph)} \left\{ \|f - g\|_p + t^r \max_{1 \le i<j \le d}
        \|D_{i,j}^r g\|_p\right\}.
\end{equation}
\end{defn}

As in the classical setting, these two gadgets are equivalent (\cite[Theorem 3.6]{DaX}).

\begin{thm} \label{thm:equivalent}
Let $r \in \NN$ and let $f \in L^p(\sph)$ if $1 \le p < \infty$
and $f \in C(\sph)$
 if $p=\infty$.  For $0 < t < 1$,
$$
   \o_r(f,t)_p \sim K_r(f,t)_p, \qquad 1 \le p \le \infty.
$$
\end{thm}

Finally, we point out that there are several  well studied moduli of smoothness on
$\sph$ (see, for instance, \cite{Di1, Rus}). One of the advantages of our new modulus
is that it reduces many problems in approximation on $\sph$ to the corresponding
problems of trigonometric approximation of one variable, the latter is classical and
well studied. Another advantage is that our modulus is relatively easier to compute,
as demonstrated in Part 3 of \cite{DaX}.

\section{ Sobolev spaces and simultaneous approximation on $\sph$}
\setcounter{equation}{0}

The classical Sobolev space $W_p^r$ on $\sph$ is defined  via the fractional
order Laplace-Beltrami operator (see, for example, \cite{BD, LN,Rus}):
\begin{equation}\label{Sobolev1}
W_p^r:= \left\{ f\in L^p(\sph):\   \|f\|_{W_p^r}:=\|f\|_p+\|(-\Delta_0)^{r/2} f\|_p<\infty\right\}.
\end{equation}
We shall introduce a new Sobolev type space on $\sph$ in this
section and then study approximation by polynomials for functions
in this new space. Our  new Sobolev space is defined via the
differential operators $D_{i,j}$, $1\leq i<j\leq d$, which, by
\eqref{Laplace-Betrami2}, are more primitive than $\Delta_0$.
First, however, we need a lemma.

\begin{lem} \label{lem:parts-Din}
For $f, g \in C^1(\sph)$ and $1\leq i\neq j\leq d$,
\begin{equation} \label{eq:parts-Dij}
   \int_{\sph} f(x) D_{i,j} g(x) d\s(x) =  -  \int_{\sph} D_{i,j} f (x)g(x)d\s(x).
\end{equation}
\end{lem}

\begin{proof}
By the rotation invariance of the Lebesgue measure $d\s$, we obtain, for
any $\t\in [-\pi,\pi]$,
$$
  \int_{\sph} f(x) g (Q_{i,j,-\t} x)\, d\s(x)=\int_{\sph}f(Q_{i,j,\t} x) g(x)\, d\s(x).
$$
Differentiating  both sides of this identity with respect to $\t$ and evaluating the
resulted equation at $\t=0$ lead to the desired equation \eqref{eq:parts-Dij}.
\end{proof}

The equation \eqref{eq:parts-Dij} allows us to define distributional derivatives
$D_{i,j}^r$ on $\sph$ for $r\in\NN$ via the identity,
\begin{equation*}
   \int_{\sph} D_{i,j}^r f(x) g(x) d\s(x) =  (-1)^r  \int_{\sph}  f(x)D_{i,j}^r
   g(x)d\s(x),\   \   g\in C^\infty(\sph).
\end{equation*}
We can now define our new Sobolev space on the sphere.

\begin{defn}
For $r\in\NN$ and $1\leq p\leq \infty$, we  define the Sobolev space
$\CW_p^r \equiv \CW_p^r(\sph)$ to be the space of functions $f\in L^p(\sph)$
whose distributional derivatives $D_{i,j}^r f$, $1\leq i<j\leq d$, all belong to
$L^p(\sph)$, with norm
$$
   \|f\|_{\mathcal{W}_p^r(\sph)} :=\|f\|_p+\sum_{1\leq i<j\leq d}\|D_{i,j}^r f\|_p,
$$
where $L^p(\sph)$ is replaced by $C(\sph)$ when $p=\infty$.
\end{defn}

The following proposition compares the new Sobolev space with the classical
one defined in \eqref{Sobolev1}:

\begin{prop}\label{prop-2-3}
For $1<p<\infty$ and $r=1$ or $2$, one has
$$
   W_p^r =\CW_p^r \quad  \text{and} \quad  \|f\|_{\CW_p^r}\sim \|f\|_{W_p^r}.
$$
In general,  for $r\ge 3$ and $1<p<\infty$,
$$
W_p^r\subset \CW_p^r \quad   \text{and}\quad   \|f\|_{\CW_p^r} \ll  \|f\|_{W_p^r}.
$$
\end{prop}

\begin{proof}
These are immediate consequences of (3.13) and (3.17) of \cite{DaX}.
\end{proof}

\begin{thm}
If  $r \in \NN$,   $f \in \CW_p^{r}$, and $1 \le p \le \infty$, then
\begin{equation}\label{eq:best-Dij}
E_{2n}(f)_p \le c n^{-r} \max_{1 \le i < j \le d} E_n(D_{i,j}^r f)_p.
\end{equation}
Furthermore,  $V_n f$, defined by \eqref{Vnf}, provides the near best
simultaneous approximation for all $D_{i,j}^r f$, $1\leq i<j\leq d$, in the sense that
\begin{equation} \label{eq:best-best-Dij}
   \|D_{i,j}^r (f - V_n f)\|_p \le c  E_n ( D_{i,j}^r f)_p, \quad 1 \le i < j \le d.
\end{equation}
\end{thm}

\begin{proof}
Applying \eqref{eq:parts-Dij} to the function $g: [-1,1] \mapsto \RR$ and using
$$
    D_{i,j}^{(x)} g(\la x, y \ra) = g' (\la x,y\ra) (x_i y_j - x_j y_i) = - D_{i,j}^{(y)} g(\la x , y \ra),
$$
it follows immediately that
$$
  D_{i,j}^{(x)} \int_{\sph} f(y) g( \la x, y \ra) d\s(y)=   \int_{\sph} D_{i,j} f (y)  g(\la x, y\ra)d\s(y).
$$
Consequently, by \eqref{zonal} and \eqref{Vnf}, we see that $V_n D_{i,j}^r = D_{i,j}^r V_n$.
Thus, using Theorems \ref{thm:best} and \ref{thm:equivalent}, we obtain
\begin{align*}
   E_{2n} (f)_p & = E_{2n} (f - V_n f)_p \le c K_r( f - V_n f, n^{-1})_p  \\
                   & \le c n^{-r} \max_{1 \le j < j \le d} \|D_{i,j}^r ( f - V_n f )\|_p \\
                   & = c n^{-r} \max_{1 \le j < j \le d} \|D_{i,j}^r  f - V_n (D_{i,j}^r f )\|_p \\
                   & \le c n^{-r} \max_{1 \le i < j \le d} E_n (D_{i,j}^r f)_p,
\end{align*}
where we used \eqref{eq:K-func-sphere} in the third step, the fact
that $V_n D_{i,j}^r=D_{i,j}^r V_n$ in the fourth step, and \eqref{Vnf-best} in the
last step. This proves \eqref{eq:best-Dij}. The inequality \eqref{eq:best-best-Dij} 
follows immediately from the above proof.
\end{proof}

Next we define a Lipschitz space on the sphere and consider approximation in such
a space.

\begin{defn}
For $r\in\NN$,  $1\leq p\leq \infty$, and $\a\in [0,1)$, we define
the Lipschitz space $\CW_p^{r,\a}\equiv \CW_p^{r,\a}(\sph)$ to be
the space of all functions $f\in \CW_p^r$ with
$$
    \|f\|_{\mathcal{W}_p^{r,\a}(\sph)} :=\|f\|_{p} + \max_{1\leq i<j\leq d} \sup_{0<|\t|\leq 1}
       \f{\|\Delta_{i,j,\t}^\ell(D_{i,j}^r f)\|_p}{|\t|^\a}<\infty,
$$
where $\ell$ is a fixed positive integer, for example, $\ell =1$.
\end{defn}

Our next theorem gives an equivalent characterization of the space
$\CW_p^{r,\a}$. For the same set of parameters as in the
definition of $\CW_p^{r,\a}$, we define the  space
$$
 H_p^{r+\a}:= \left \{ f\in L^p(\sph):  \|f\|_{H_p^{r+\a}}: = \|f\|_p +
      \sup_{0<t\leq 1} \f{\o_{r+\ell}(f, t)_p} {t^{r+\a}}<\infty \right\}.
$$

\begin{thm}\label{thm-2-5}
If $r\in\NN$,  $1\leq p\leq \infty$, and $\a\in [0,1)$, then
$$
  \CW_p^{r,\a} = H_p^{r+\a} \quad \hbox{and}\quad  \|f\|_{\CW_p^{r,\a}} \sim \|f\|_{H_p^{\a+r}}.
$$
\end{thm}

\begin{proof}
To prove that $f \in \CW_p^{r,\a}$ implies $f\in H_p^{r+\a}$ and
$\|f\|_{H_p^{r+\a}} \leq c \|f\|_{\CW_p^{r,\a}}$, it suffices to show
that for $f\in \CW_p^r$ and $\ell\in \NN$,
\begin{equation}\label{2-15-1}
    \|\tr_{i,j,\t}^{r+\ell} f \|_p \leq c |\t|^r\|\tr_{i,j, \t}^\ell (D_{i,j}^r f)\|_p.
\end{equation}
Using Lemma 2.6 (ii) in \cite{DaX}, we have
\begin{align*}
    \|\tr_{i,j,\t}^{r+\ell} f \|_p =\left \|\tr_{i,j,\t} ^r (\tr_{i,j,\t}^\ell f)
          \right \|_p\leq c  |\t|^r \left \|D_{i,j}^r (\tr_{i,j,\t}^\ell f)\right \|_p.
\end{align*}
However, from \eqref{Delta_ij} and \eqref{partial_ij2}, a quick computation
shows that $\tr_{i,j, \t} D_{i,j} =D_{i,j} \tr_{i,j,\t}$, hence, by iteration,
$\tr_{i,j, \t} ^\ell D_{i,j}^r =D_{i,j}^r \tr_{i,j,\t}^\ell$. As a result,
$$
\|D_{i,j}^r \tr_{i,j,\t}^\ell f \|_p= \|\tr_{i,j,\t}^\ell D_{i,j}^r f\|_p.
$$
Together, these two displayed equations yield \eqref{2-15-1}.

Conversely, assume $f\in H_p^{r+\a}$. We first show that $D_{i,j}^r f \in L^p(\sph)$.
For $g\in C^\infty(\sph)$ and $G_x (t) := g(Q_{i,j,t} x)$, we can write \cite[(4.8)]{DaX}
$$
\tr_{i,j,\t}^r (g)(x) =\int_0^\t \cdots\int_0^\t G_x^{(r)}(t_1+\cdots+t_r)\, dt_1\cdots dt_r,
$$
which implies, in particular, that
\begin{equation}\label{2-16}
  \lim_{\t\to 0} \f{\tr_{i,j,\t}^r (g)(x)}{\t^r}=G_x^{(r)}(0)=D_{i,j}^r g(x).
\end{equation}
Thus, by the definition of the distributional derivative $D_{i,j}^r f$, it follow
that, for $g\in C^\infty(\sph)$,
\begin{align*}
&\int_{\sph} [D_{i,j}^r f(x)] g(x)\, d\s =(-1)^r \int_{\sph} f(x)D_{i,j}^r g(x) \, d\s \\
&  = (-1)^r \lim_{\t\to 0} \int_{\sph} f(x) \f {\tr_{i,j,\t}^r g (x)}{\t^r}\, d\s
   =(-1)^r \lim_{\t\to 0} \int_{\sph} \f {\tr_{i,j,-\t}^r f(x)}{\t^r} g(x)\, d\s,
\end{align*}
where the last step uses the rotation invariance of the Lebesgue measure $d\s$.
However, by the Marchaud inequality (Proposition 2.7 of \cite{DaX}),
for $f \in  H_p^{r+\a}$,
\begin{align*}
  \o_{r} (f, t)_p \leq   c_\ell  t^r \int_t^1 \frac{\o_{r + \ell} (f, u)_p}{u^{r+1}} du
  \le c t^r \|f\|_{H_p^{r+\a}}.
\end{align*}
Hence, by H\"older's inequality, we deduce with $\frac{1}{p} + \frac{1}{p'} =1$ that
$$
\left | \int_{\sph}[  D_{i,j}^r f(x)] g(x)\, d\s(x)\right | \leq c \|f\|_{H_p^{r+\a}} \|g\|_{p'},
$$
which implies, upon taking supreme over all $g$ with $\|g\|_{p'} \le 1$ that
$D_{i,j}^r f\in L^p(\sph)$. Next we note that for $\ell\in \NN$,
\begin{equation}\label{2-17}
 \|\tr_{i,j,\t}^\ell (D_{i,j}^rf)\|_p \leq c \int_0^\t \|\tr_{i,j,u}^{r+\ell} f\|_p \f {du}{u^{r+1}},
\end{equation}
which follows from the analogue result for trigonometric functions \cite[(7.1)]{Di0}
as in the proof of Lemma 2.6 of \cite{DaX}. Consequently, it follows that
$$
   \|f\|_{\CW_p^{r,\a}} \leq \|f\|_p + c \max_{1 \le i < j\le d} \sup_{0 \le |\t| \le 1}
         \frac{1}{|\t|^\a} \int_0^\t \|\tr_{i,j,u}^{r+\ell} f\|_p \f {du}{u^{r+1}}
        \le c  \|f\|_{H_p^{r+\a}}
$$
since $0 < \a < 1$. This completes the proof.
\end{proof}

Theorem \ref{thm-2-5} together with the Jackson theorem implies the
following:

\begin{cor}
If  $r \in \NN$,  $\a\in [0,1)$,   $f \in \CW_p^{r,\a}$, and $1
\le p \le \infty$, then
$$
   E_n(f)_p \leq  c  n^{-r-\a}\|f\|_{\CW_p^{r,\a}}.
$$
\end{cor}

\section{Polynomial approximation on $\BB^d$: Recent Progress}
\setcounter{equation}{0}

In this section, we recall recent progress on polynomial approximation on
$\BB^d$ as developed in \cite{DaX}. For $\mu\ge 0$, let $W_\mu$ denote
the weight function  on $\BB^d$ defined by
\begin{equation}\label{weightW}
   W_\mu(x) : = (1 - \|x\|^2)^{\mu -1/2}.
\end{equation}
For $1 \le p < \infty$ we denote by $\|f\|_{p,\mu}$ the norm for the weighted
$L^p$ space  $L^p(\BB^d, W_\mu)$,
\begin{equation}
  \|f\|_{p,\mu}: = \left( \int_{\BB^d} |f(x)|^p W_\mu(x) dx \right)^{1/p},
\end{equation}
and $\| f\|_{\infty,\mu} : = \|f\|_\infty$ for $f \in C(\BB^d)$.
When we need to emphasis that the norm is taken over $\BB^d$, we
write $\|f\|_{p,\mu} = \|f\|_{L^p(\BB^d,W_\mu)}$. For $f \in L^p(\BB^d, W_\mu)$,
$1 \le p< \infty$, or $f \in C(\BB^d)$, $p = \infty$, the best approximation by
polynomials is defined by
$$
       E_n(f)_{p,\mu}: = \inf_{p \in \Pi_n^d} \|f - p \|_{p,\mu}.
$$

\subsection{Weighted orthogonal polynomial expansions on $\BB^d$}
Let $\CV_n^d(W_\mu)$ denote the space of orthogonal polynomials of
degree $n$ with respect to the weight function $W_\mu$ on $\BB^d$.
We denote by $P_n^\mu(x,y)$ the reproducing kernel of $\CV_n^d(W_\mu)$
in $L^2(\BB^d, W_\mu)$. It is shown in \cite[Theorem 2.6]{X01} that
\begin{equation}\label{reprodB-S}
 P_n^\mu (x,y) = \int_{\SS^{m-1}} Z_{n,d+m} \left (  \la x, y \ra
            +\sqrt{1-\|y\|^2} \, \la x', \xi \rangle \right ) d\s(\xi)
\end{equation}
for any $x, y \in \BB^d$ and  $(x,x') \in \SS^{d+m-1}$,  where $Z_{n,d}(t)$ is
the zonal harmonic defined in \eqref{zonal} and $\mu = \frac{m-1}{2}$.
For $\eta$ being a $C^\infty$-function on
$[0,\infty)$ that satisfies the  properties as defined in Section 1.1, we define
an operator
\begin{equation}\label{Vnf_mu}
 V_n^\mu f(x) := a_\mu \int_{\BB^d} f(y) K_n^\mu(x,y) W_\mu(y) dy, \quad x \in \BB^d,
\end{equation}
where $a_\mu$ is the normalization constant of $W_\mu$ and
\begin{equation}\label{kernelB}
   K_n^\mu(x,y): = \sum_{k=0}^{2n} \eta \left(\frac{k}{n}\right) P_k^\mu(x,y).
\end{equation}
This operator plays the same role as $V_n f$ in the study on $\sph$. In
particular, $V_n^\mu f \in \Pi_{2n}^d$ and $\|f - V_n^\mu f \|_{p,\mu} \le c
E_n (f)_{p,\mu}$, $1 \le p \le \infty$.

The spaces $\CV_n^d(W_\mu)$ of orthogonal polynomials  are also
the eigenspaces of the following second order differential
operator:
\begin{equation}\label{D-mu}
  \CD_\mu:=  \sum_{i=1}^d (1-x_i^2)  \partial^2_i - 2
   \sum_{1\le i < j \le d}  x_i x_j \partial_i \partial_j  -
     (d+2 \mu)\sum_{i=1}^d x_i \partial_i.
\end{equation}
Indeed, elements of $\CV_n^d(W_\mu)$ satisfy (cf. \cite[p. 38]{DX})
\begin{equation} \label{D-eigen}
\CD_\mu P =  - n(n+ d+ 2\mu-1 ) P \quad\hbox{for all} \quad P \in
\CV_n^d(W_\mu).
\end{equation}
It was shown in \cite[Proposition 7.1]{DaX}  that the differential
operator $D_\mu$ can be decomposed as a sum of second order
differential operators:
\begin{equation}\label{decomp}
  \CD_\mu = \sum_{i=1}^d D^2_{i,i} +  \sum_{1\le i < j \le d} D^2_{i,j}
      =  \sum_{1\le i \le j \le d} D^2_{i,j},
\end{equation}
where the operators $D^2_{i,j}$ for $1 \le i < j \le d$  are defined as
in \eqref{partial_ij}, and
\begin{equation}
   D^2_{i,i} \equiv D_{\mu,i,i}^2:= [W_\mu(x)]^{-1} \partial_i \left[(1-\|x\|^2) W_\mu(x) \right]
        \partial_i, \qquad 1 \le i \le d.
\end{equation}
%It can be easily verified that the operators $D_{i,j}^2$, $1\leq i, j\leq d$ commute
%with $\CD_\mu$, which, in particular, implies that the spaces $\CV_n^d(W_\mu)$
%of orthogonal polynomials  are invariant under the operators $D_{i,j}$, $1\leq i\leq j
%\leq d$.

For the rest of the paper, we will always set $\varphi (x)=
\sqrt{1-\|x\|^2}$. We have the following useful estimates:
\begin{prop}
\label{prop-3-1}If  $1<p<\infty$,  and $g\in C^2(\BB^d)$ then
\begin{equation}\label{3-27}
   \|\CD_\mu  g\|_{p,\mu} \sim \sum_{1 \le i \le j \le d} \|D^2_{i,j}g
   \|_{p,\mu}.
\end{equation}
Furthermore, if $r\in \NN$, $1<p<\infty$, and $g\in C^{2r}(\BB^d)$
then
\begin{equation}\label{3-28}
  c_1 \| \varphi^{2r} \partial_{i}^{2r} g\|_{p,\mu} \le \|D^{2r}_{i,i} g\|_{p,\mu}
      \le  c_2 \| \varphi^{2r} \partial_{i}^{2r} g\|_{p,\mu} + c_2
      \|g\|_{p,\mu}, \   \  1\leq i\leq d.
\end{equation}\end{prop}

\begin{proof}\eqref{3-27} was proved in \cite[Theorem 7.3]{DaX}. In
the case when $r=1$, \eqref{3-28} was shown in \cite[Theorem
7.4]{DaX}, and the proof there works equally well for $r\ge 1$.
\end{proof}

\subsection{Moduli of smoothness and K-functionals}
Two moduli of smoothness and their equivalent K-functionals on $\BB^d$
were introduced in \cite{DaX}.  Our first modulus on $\BB^d$ was defined via
an extension $\wt{f}$  of a function $f:\BB^d\to\RR$ to $\BB^{d+1}$, defined by
\begin{equation} \label{wt-f}
  \wt f(x,x_{d+1}) : = f(x), \qquad (x,x_{d+1}) \in \BB^{d+1}, \quad x \in \BB^d.
\end{equation}
More precisely, it is defined as follows (\cite[Definition 5.3]{DaX}).

\begin{defn}
For $r \in \NN$, $t > 0$, and $f \in L^p(\BB^d, W_\mu)$, $1 \le p < \infty$, or
$f \in C(\BB^d)$ for $p = \infty$, define
\begin{align}  \label{omegaB-def}
  \o_r(f,t)_{p,\mu} :=  \sup_{|\t|\le t} &
  \left\{ \max_{1 \le i<j \le d}  \| \tr_{i,j,\t}^r  f\|_{L^p(\BB^{d},W_ {\mu})}, \right. \\
& \quad
 \left. \max_{1 \le i \le d}  \| \tr_{i,d+1,\t}^r \wt f\|_{L^p(\BB^{d+1},W_ {\mu-1/2})} \right\},
\notag
\end{align}
where   for $m=1$, $ \| \tr_{i,d+1,\t}^r \wt f\|_{L^p(\BB^{d+1},W_
{\mu-1/2})}$ is replaced by  $\| \tr_{i,d+1,\t}^r \wt f\|_{L^p(\SS^d)}$.
\end{defn}

In the case when $\mu = \frac{m-1}{2}$ and $m\in \NN$,  we have established
in \cite[Theorem 5.5]{DaX} the direct theorem, that is, the Jackson inequality
\begin{equation}\label{JacksonB}
   E_n (f)_{p,\mu} \le c\, \o_r(f,n^{-1})_{p,\mu}, \qquad 1 \le  p \le  \infty,
\end{equation}
and the corresponding inverse theorem,
\begin{equation} \label{inverseB}
\o_r (f,n^{-1})_{p,\mu} \le c\, n^{-r} \sum_{k=1}^n k^{r-1}E_{k}(f)_{p,\mu},
\end{equation}
in terms of this new modulus of smoothness. Moreover, it was also shown in
\cite[Theorem 5.8]{DaX} that the modulus of smoothness $\o_r (f,t)_{p,\mu}$  is
equivalent to the following K-functional:
\begin{align}  \label{K-functB-def}
 &   K_r(f,t)_{p,\mu}
:= \inf_{g\in C^r(\BB^d)} \Big\{  \|f-g\|_{L^p(\BB^d, W_\mu)}   \\
  & \qquad
     + t^r   \max_{1 \le i<j \le d}  \| D_{i,j}^r  g\|_{L^p(\BB^{d},W_ {\mu})} +
   t^r \max_{1 \le i \le d}  \| D_{i,d+1}^r  \wt g\|_{L^p(\BB^{d+1},W_ {\mu-1/2})} \Big \}.
\notag
\end{align}
where if $m =1$, then $\| D_{i,d+1}^r  \wt g\|_{L^p(\BB^{d+1},W_{\mu-1/2})}$ is
replaced by $\| D_{i,d+1}^r  \wt g\|_{L^p(\SS^d)}$.

Our second new modulus  in \cite{DaX} can be considered as a
higher-dimensional analogue of  the classical Ditzian-Totik modulus on the
interval $[-1,1]$. In the unweighted case, this modulus is defined as follows
\cite[Definition 6.7]{DaX}:

\begin{defn} \label{defn:modulusB2}
For $r \in \NN$, $t > 0$, and $f \in L^p(\BB^d, W_\mu)$, $1 \le p < \infty$, or
$f \in C(\BB^d)$ for $p = \infty$, define
\begin{equation}\label{modulusB2}
  \wh \o_r(f,t)_{p} := \sup_{0 < |h| \le t} \left\{ \max_{1 \le i<j \le d}
      \| \tr_{i,j,h}^r  f\|_{p},
      \max_{1 \le i \le d}  \|\wh \tr_{h\varphi  e_i}^r f\|_{p} \right
      \},
\end{equation}
where $\wh\Delta_h$ denotes the central difference and
$$
\wh\Delta_{h \varphi e_i}^r f (x) : = \sum_{k=0}^r (-1)^k
\binom{r}{k} f\left(x+ (\tfrac{r}2 -k)h\varphi (x) e_i\right),
$$
and we assume $\wh\Delta_{h\varphi  e_i}^r = 0$ if either
of the points $x \pm r \frac{h\varphi (x)}{2} e_i$ does not belong
to $\BB^d$.
\end{defn}

The weighted version  $\wh \o_r(f,t)_{p,\mu}$ of the above modulus
can also be defined, but is more complicated. Both direct and inverse
theorems were established in terms of $\wh \o_r(f,t)_{p}$ in
\cite[Theorem 6.13]{DaX}. Furthermore, it was shown in
\cite[Theorem 6.10]{DaX} that the modulus of smoothness
$\wh \o_r (f,t)_{p}$  is equivalent to the following K-functional
\begin{align*}  %\label{K-functB2-def}
  \wh K_r(f,t)_{p,\mu}:= \inf_{g\in C^r(\BB^d)} \Big\{ \|f-g\|_{p,\mu}
        + t^r   \max_{1 \le i<j \le d}  \| D_{i,j}^r  g\|_{p,\mu} +
          t^r \max_{1 \le i \le d}  \| \varphi^r \partial_i^r g\|_{p,\mu} \Big \}
\end{align*}
in the sense that, for the equivalence between $\wh K_r(f,t)_{p}$ and $ \o_r(f,t)_{p}$
in the unweighted case,
$$
c^{-1}\wh\o^r(f,t)_p\le  \wh K_r(f,t)_p\leq c\, \wh\o^r(f,t)_p+c\,t^r\|f\|_p.
$$

The two $K$-functionals, hence their equivalent moduli of smoothness,
are connected as shown in \cite[Theorem 6.2]{DaX}.

\begin{thm}\label{thm:whK=K}
Let $\mu = \frac{m-1}{2}$ and $m \in \NN$. Let $f \in
L^p(\BB^d,W_\mu)$ if $1 \le p < \infty$, and $f\in C(\BB^d)$ if
$p=\infty$. We further assume that $r$ is odd when $p=\infty$.
Then
\begin{align} \label{whK=K_r=1}
   \wh  K_1(f,t)_{p,\mu} \sim K_1(f,t)_{p,\mu},
\end{align}
and for $r \ge  1$, there is a $t_r > 0$ such that
\begin{align} \label{whK=K}
       K_r(f,t)_{p,\mu} \le c \, \wh K_r(f,t)_{p,\mu}
             + c \, t^r \|f\|_{p,\mu},  \qquad 0<t<t_r.
\end{align}
\end{thm}

Finally, we point out that it was shown in \cite{DaX} that both moduli
$\wh \o_r(f,t)_{p,\mu}$ and $ \o_r(f,t)_{p,\mu}$ enjoy most of the properties
of classical moduli of smoothness and they are computable as demonstrated
in Part 3 of \cite{DaX}. In comparison, the only other modulus of smoothness
\cite{X05} on the unit ball that is strong enough to characterize the best approximation
is hardly computable.

\subsection{Representation of the term  $D_{i,d+1}^r\wt{f}$}
The term  $D_{i,d+1}^r \wt g$ appears in the definition of our first
$K$-functional $K_r(f,t)_{p,\mu}$ in \eqref{K-functB-def} on the ball,
where $\wt g(x,x_{d+1}) = g(x)$ as in \eqref{wt-f}. Notice that $\wt f$
is a function in $x \in \RR^d$, but the operator $D_{i,d+1} = x_i \partial_{d+1}
 - x_{d+1} \partial_i$ involves $x_{d+1}$, so that $D_{i,d+1}^r \wt f$ is
indeed a function of $(x,x_{d+1})$ in $\BB^{d+1}$. The following lemma
gives an explicit formula of  this term in terms of $f$.

\begin{lem} \label{lem:D_i_d+1}
Assume that  $(y, y_{d+1})=s(x, x_{d+1}) \in \BB^{d+1}$ with $s=\|(y, y_{d+1})\|>0$,
$x\in \BB^d$ and $x_{d+1}=\vi(x)\ge 0$. If $f\in C^r(\BB^d)$, then
$$
(D_{i, d+1}^r \wt f ) (y, y_{d+1})  =  \left( - \varphi (x) \frac{\partial}{\partial x_i} \right)^r \Bigl[f(s x)\Bigr],
     \quad 1 \le i \le d.
$$
\end{lem}

\begin{proof}
The proof uses induction. For $r = 1$, we have
$$
  D_{i,d+1} \wt f(y, y_{d+1}) = (y_i \partial_{d+1} - y_{d+1} \partial_i) f(y) = - y_{d+1} \partial_i f(y).
$$
Hence, using the fact that $\frac{\partial}{\partial x_i} [f(s x)] =
s (\partial_i f)(s x)$ we have
$$
  (D_{i,d+1} \wt f) ( s x, s x_{d+1}) =  - s x_{d+1} (\partial_i f)(sx) = -s \varphi(x) (\partial_i f)(s x)
     =  - \varphi (x) \frac{\partial}{\partial x_i} [f(s x)].
$$
Let $F_r(x,x_{d+1})= D_{i, d+1}^r \wt f(x,x_{d+1})$. Assume that the
result has been established for $r$. Then $F_r (s x,  s \varphi(x))
= (-\varphi \partial_{i})^r [f(s x)]$. By definition,
\begin{align} \label{eq:F_r+1}
 F_{r+1}(s x, s x_{d+1}) & = (D_{i,d+1} F_r)(s x,s x_{d+1})  \\
   &  = s x_i (\partial_{d+1} F_r)(s x,s x_{d+1}) - s x_{d+1} (\partial_i F_r)( s x, sx_{d+1}). \notag
\end{align}
On the other hand, taking derivative by chain rule shows that
\begin{align*}
 \left( - \varphi (x) \frac{\partial}{\partial x_i} \right)^{r+1} [f(s x)] & =
     \left( - \varphi (x) \frac{\partial}{\partial x_i} \right) [F_r(s x, s \varphi(x))]\\
    &  = - s \varphi (x) (\partial_i F_r)(s x, s \varphi(x)) +  s x_i (\partial_{d+1} F_r)(s x, s \varphi(x)),
\end{align*}
which is the same as the right hand side of \eqref{eq:F_r+1} with
$x_{d+1} = \varphi(x)$.
\end{proof}

\begin{lem} \label{lem:parity_Fr}
The function $D_{i,d+1}^r \wt f(x,x_{d+1})$ is even in $x_{d+1}$ if
$r$ is even, and odd in $x_{d+1}$ if $r$ is odd.
\end{lem}

\begin{proof}
For $r = 1$, $D_{i,d+1} \wt f(x,x_{d+1}) = - x_{d+1} \p_if(x)$ is
clearly odd in $x_{d+1}$. And
$$
   D_{i,d+1}^2 \wt f(x,x_{d+1}) = - x_i \partial_i f(x) + x_{d+1}^2 \partial_i^2 f(x)
$$
is even in $x_{d+1}$. The general case follows from induction upon
using \eqref{eq:F_r+1}.
\end{proof}

Recall that our $K$-functional $K_r(f,t)_{p,\mu}$ with $\mu  = \frac{m-1}{2}$
in \eqref{K-functB-def} is defined. when $m =1$, with $\|D^r_{i,d+1} \wt g\|_{L^p(\SS^d)}$ in
place of $\|D^r_{i,d+1} \wt g\|_{L^p(\BB^{d+1}, W_{\mu-1/2})}$. Hence, as a
consequence of the above lemmas, we conclude the following: 

\begin{prop}
For $g\in C^r (\BB^d)$ and  the Chebyshev weight $W_0$ on $\BB^d$,
we have
\begin{equation} \label{D_id+1W0}
   \|D_{i, d+1}^r \wt {g}\|_{L^p(\SS^{d})}=  \| (\varphi \partial_i)^r g\|_{L^p(\BB^d,W_0) }.
\end{equation}
\end{prop}

\begin{proof}
Let $\SS_+^d = \{x \in \SS^d: x_{d+1} \ge 0\}$. By Lemma
\ref{lem:parity_Fr} we only need to consider $\SS_+^d$ when dealing
with $D_{i,d+1}^r \wt g$. By Lemma \ref{lem:D_i_d+1} with $s= 1$, we
then obtain
\begin{align*}
  \int_{\SS^d} \left |D_{i,d+1}^r \wt g (x,x_{d+1})\right|^p d \s(x,x_{d+1}) &
     =  2 \int_{\SS^d_+} \left|(\varphi(x) \partial_i)^r g(x)\right |^p d\s(x,x_{d+1}) \\
     &  =  \int_{\BB^d} \left |(\varphi(x) \partial_i)^r g(x)\right |^p \frac{dx}{\sqrt{1-\|x\|^2}},
\end{align*}
which is what we want to prove.
\end{proof}

In general, using polar coordinates, Lemma \ref{lem:D_i_d+1}, and
Lemma \ref{lem:parity_Fr}, we can deduce
\begin{prop}
If $g\in C^r(\BB^d)$, $\mu=\f {m-1}2$ and $m>1$,   then for $1\leq p<\infty$,
\begin{align} \label{D_i_d+1norm}
 &  \|D_{i,d+1}^r \wt g\|_{L^p(\BB^{d+1}, W_{\mu -1/2})}^p \\
 & \qquad  =  \int_0^1 s^d (1-s^2)^{\mu-1} \int_{\BB^d}  \bigl|(\varphi(x)\partial_i )^r [g(s x)]\bigr|^p
     \frac{dx}{\sqrt{1-\|x\|^2}} ds;\notag
\end{align}
whereas for $p=\infty$, we have
$$
   \max_{y \in \BB^{d+1}} | D_{i,d+1}^r \wt g(y)|  = \max_{x \in \BB^d,0\le s \le 1}
           \left| \left(\varphi(x) \frac{\partial}{\partial x_i} \right)^r [g(s x)] \right|.
$$
\end{prop}

\section{ Sobolev Spaces and Simultaneous Approximation on $\BB^d$}
\setcounter{equation}{0}

We start with the definition of a Sobolev space on $\BB^d$.

\begin{defn}
For $1 \le p \le \infty$, $f\in C^r(\BB^d)$,  and $r\in\NN$,   we define
\begin{equation}\label{6-35}
  \|f\|_{\CW_p^r(\BB^d, W_\mu)}:=\|f\|_{p,\mu} + \sum_{1\leq i<j\leq d} \|D_{i,j}^r f\|_{p,\mu}
     + \sum_{i=1}^d\|\varphi^r \partial_i^r f\|_{p,\mu},
\end{equation}
and  define $\CW_p^r(\BB^d, W_\mu)$  to be the completion of $C^r(\BB^d)$
with respect to the norm $\|\cdot\|_{\CW_p^r(\BB^d,W_\mu)}$.
\end{defn}

\begin{rem}Since convergence in the norm $\|\cdot\|_{\CW_p^r(\BB^d,
W_\mu)}$ implies convergence in the weighted $L^p$-norm
$\|\cdot\|_{p,\mu}$, we may assume that $\CW_p^r(\BB^d,
W_\mu)\subset L^p(\BB^d, W_\mu)$ when $p<\infty$, and
$\CW_p^r(\BB^d, W_\mu)\subset C(\BB^d)$ when $p=\infty$. As a
consequence, we can also extend the definitions of the  operators
$D_{i,j}^r f$ and $\varphi^r \partial_i^r f$ to the whole space
$\CW_p^r(\BB^d, W_\mu)$.
\end{rem}

The following proposition follows readily from (6.15) and (6.16) of \cite{DaX}
and Proposition \ref{prop-3-1}:

\begin{prop}
If $f\in C^r(\BB^d)$, $\mu\ge 0$ and $1\leq p\leq \infty$, then
\begin{equation}\label{6-36}
\sum_{1\leq i \leq d} \|D_{i, d+1}^r \wt{f}\|_{L^p(\BB^{d+1}, W_{\mu-\frac12})}
        \leq c \|f\|_{\CW_p^r(\BB^d, W_\mu)}.
\end{equation}
Furthermore, if $f\in C^{2r}(\BB^d)$ and $1<p<\infty$ then
\begin{equation}
\|f\|_{\CW_p^{2r}(\BB^d, W_\mu)} \sim \sum_{1\leq i\leq j\leq d}
\|D_{i,j}^{2r} f\|_{p,\mu}.
\end{equation}
\end{prop}

\begin{thm} \label{thm-6-2}Let $\mu = \frac{m-1}{2}$ with $m\in\NN$. For $f \in \CW_p^r (\BB^d,
W_\mu)$, $1 \le p \le \infty$,
\begin{align} \label{Ball-bast-best}
   E_{2n} (f)_{p,\mu} \le c n^{-r} & \left[  \max_{1 \le i < j \le d}  E_n (D_{i,j}^r f)_{p,\mu} \right. \\
          & + \left.  \max_{1 \le i \le d} E_n ( D_{i,d+1}^r  \wt f)_{L^p(\BB^{d+1}, W_{\mu -\frac12})} \right].
   \notag
\end{align}
Furthermore,  $V_n^\mu f$, defined by \eqref{Vnf_mu}, provides the near best
simultaneous approximation for all $D_{i,j}^r f$, $1\leq i<j\leq d+1$ in the sense that
\begin{align*}
 \|D_{i,j}^r (f - V_n^\mu f)\|_{p,\mu} & \le c  E_n ( D_{i,j}^r f)_{p,\mu} \quad 1 \le i < j \le d \\
 \|D_{i,d+1}^r (\wt f - \wt{V_n^\mu f} )\|_{L^p(\BB^{d+1}, W_{\mu -\frac12})}
     & \le c  E_n( D_{i,d+1}^r  \wt f)_{L^p(\BB^{d+1}, W_{\mu -\frac12})}, \quad 1 \le i \le d.
\end{align*}
\end{thm}

\begin{proof}
For $f$ defined on $\BB^d$, we define $F(x, x'): = f(x)$, $x \in \BB^d$, $(x,x') \in \SS^{d+m-1}$.
By \cite[Lemma 5.2]{DaX}, $(V_n F) (x, x')  =  V_n^\mu f(x)$. Furthermore, by
\cite[Lemma 5.7]{DaX}, $K_r( f, n^{-1} )_{L^p(\BB^d, W_\mu)} \sim K_r( F, n^{-1})_{L^p(\SS^{d+m-1})}$.
Hence, it follows that
\begin{align*}
  E_{2n}(f)_{L^p(\BB^d)} &  \le c K_r( f -V_n^\mu f, n^{-1} )_{L^p(\BB^d, W_\mu)} \\
      &  \le c K_r( F -V_n F, n^{-1})_{L^p(\SS^{d+m-1})} \\
      &  \le c n^{-r}   \max_{1 \le i < j \le d+m} \| D_{i, j}^r (F - V_n F)\|_{L^p(\SS^{d+m-1})}\\
      & =  c n^{-r}   \max_{1 \le i < j \le d+1} \| D_{i, j}^r F -  (D_{i,j}^r V_n F\|_{L^p(\SS^{d+m-1})},
\end{align*}
where the last step follows from the fact that $V_n(D_{i,d+k}^r F)$ depends on
$x_j$, $1 \le j \le d$, and $x_{d+k}$, which implies that we only need to consider
$1 \le i < j \le d+1$.

Denote by $V_{n,d}^\mu$ the operator \eqref{Vnf_mu} associated with $W_\mu$ on $\BB^d$ and
$\wt f(x,x_{d+1}) = f(x)$. By Lemma 5.2 of \cite{DaX},
\begin{equation} \label{VnF-Vnf}
  V_n F (x,x') = V_{n,d+1}^{\mu-1/2} \wt f (x,x_{d+1})= V_{n,d}^\mu f(x).
\end{equation}
Since $D_{i,j}^r V_n  = V_n D_{i,j}^r$ on the sphere, it follows that, for $1 \le i, j \le d$,
$$
 D_{i,j}^r V_{n,d}^\mu f(x) = D_{i,j}^r (V_n F) (x,x') = V_n D_{i,j}^r F(x,x') = V_{n,d}^{\mu} D_{i,j}^r f (x).
$$
Consequently, it follows from \cite[(5.8)]{DaX} that
\begin{align*}
   \|D_{i,j}^r F -  D_{i,j}^r V_n F) \|_{L^p(\SS^{d+m-1})} &
      = c  \| D_{i,j}^r f-  D_{i,j}^r V_{n,d}^\mu  f ) (x)\|_{L^p(\BB^d, W_{\mu})}\\
        & = c  \| D_{i,j}^r f- V_{n,d}^\mu (D_{i,j}^r f ) (x)\|_{L^p(\BB^d, W_{\mu})}\\
        & \le c E_n( D_{i,j}^r f)_{L^p(\BB^d, W_\mu)}.
\end{align*}
Whereas for $D_{i,d+1}^r F$ term, we have for $1 \le i \le d$,
\begin{align*}
   D_{i,d+1}^r V_{n,d+1}^{\mu -1/2} \wt f(x,x_{d+1}) & = D_{i,d+1}^r (V_n F) (x,x')\\
      & = V_n (D_{i,d+1}^r F) (x,x') =
      V_{n,d+1}^{\mu-1/2} D_{i,j}^r \wt f (x,x_{d+1}).
\end{align*}
Consequently,
\begin{align*}
    \|D_{i,d+1}^r F- D_{i,d+1}^r V_n F \|_{L^p(\SS^{d+m-1})}  & = c \| D_{i,j}^r f -
        V_{n,d+1}^{\mu-1/2} D_{i,j}^r f (x)\|_{L^p(\BB^{d+1}, W_{\mu-1/2})} \\
       & \le c E_n( D_{i,j}^r \wt f)_{L^p(\BB^{d+1}, W_{\mu-1/2})}.
 \end{align*}
This proves \eqref{Ball-bast-best}. The conclusion that $V_n^\mu f$ is the near best
simultaneous approximation follows from the above proof and \eqref{VnF-Vnf}.
\end{proof}

\begin{cor}\label{cor-6-3}
 Let $\mu = \frac{m-1}{2}$ with $m\in\NN$. If  $f \in \CW_p^r
(\BB^d, W_\mu)$, $1 \le p \le \infty$, then
$$
   E_{n} (f)_{p,\mu} \le c n^{-r}\|f\|_{\CW_p^r
(\BB^d, W_\mu)}.$$\end{cor}

\begin{proof}
This follows immediately from \eqref{6-35}, \eqref{6-36} and Theorem \ref{thm-6-2}.
\end{proof}

In the next corollary, we replace $D_{i,d+1}^r \wt f$ term in \eqref{Ball-bast-best}
by ordinary derivatives of $f$. First we consider the Chebyshev weight $W_0$ on $\BB^d$
(with $\mu=0$).

\begin{cor}
For $f \in \CW_p^r(\BB^d, W_0)$, $r\in\NN$ and  $1 \le p \le \infty$,
\begin{align*}
   E_{2n} (f)_{p,0}& \le c n^{-r}   \max_{1 \le i < j \le d}  E_n (D_{i,j}^r f)_{p, 0}
   +c n^{-r}\max_{1 \le i \le d}
E_n \left( (\vi\p_i)^r f\right)_{p,0}.
\end{align*}
\end{cor}

\begin{proof}
By \eqref{D_id+1W0}, for all $f\in C^r (\BB^d)$,
$$ ( \vi(x) \p_i)^r f (x)=D_{i, d+1}^r \wt{f} (x, x_{d+1}),$$
where $x\in \BB^d$ and $x_{d+1}=\vi(x)$. The desired conclusion
then follows.
\end{proof}

For $\mu > 0$, including the case $\mu =1/2$ (the constant weight function), however,
the best that we can do is the following:

\begin{cor} \label{cor:2.3}
Let $\mu = \frac{m-1}{2}$ and $m\in\NN$. For $f \in \CW_p^r(\BB^d, W_\mu)$, $1 \le p \le \infty$,
\begin{align*}
   E_{2n} (f)_{p,\mu} \le c n^{-r}   \max_{1 \le i < j \le d}  E_n (D_{i,j}^r f)_{p, \mu}
   & +c n^{-r}\max_{1 \le i \le d}\Bigl[
          \max_{1\leq j< \f{r+1}2}E_{n-r} (\p_i^j f)_{p,\mu}\\
   & + \max_{\f{r+1}2\leq j\leq r} E_{n-r}(\p_i^j f)_{p, \mu+(j-\f r2)p}\Bigr].
\end{align*}
\end{cor}

\begin{proof}
It was shown in Lemma 6.4 of \cite{DaX} that
$$
D_{i,d+1}^r \wt f (x,x_{d+1}) = \sum_{j=1}^r p_{j,r}(x_i,x_{d+1}) \partial_i^j f(x),
     \quad x \in \BB^d, \, (x, x_{d+1}) \in \BB^{d+1},
$$
where  $p_{j,r}$ is a polynomial of degree $\leq j$.  Since $E_n (f)$ is subadditive,
it follows that
\begin{equation}\label{6-38}
E_{n}( D_{i,d+1}^r \wt f )_{L^p(\BB^{d+1}, W_{\mu-1/2})} \le  \sum_{j= 1}^r
 \inf_{g\in\Pi_{n-r}^d}\left \|p_{j,r}  (\partial_i^j \wt f-\wt{g})\right \|_{L^p(\BB^{d+1}, W_{\mu-1/2})}.
\end{equation}
However, using (6.11) and (6.12) of \cite{DaX},  we have, for $(x, x_{d+1})\in \BB^{d+1}$,
\begin{equation*}
|p_{j,r}(x, x_{d+1})|\leq \begin{cases} c, & \text{if $1\leq j < \f {r+1}2$},\\
c |x_{d+1}|^{2j-r}, &  \text{if $ \f {r+1}2\leq j \leq r$}.   \end{cases}
\end{equation*}
Thus, by \eqref{6-38}, we deduce \begin{align*}
  E_{n}( D_{i,d+1}^r \wt f )_{L^p(\BB^{d+1}, W_{\mu-1/2})} & \leq
    c \max_{1\leq j< \f{r+1}2}\inf_{g\in \Pi_{n-r}^d} \|\p_i^j f- g\|_{p,\mu}\\
&  +c \max_{\f{r+1}2\leq j\leq r}\inf_{g\in \Pi_{n-r}^d}
\|\p_i^j f-g\|_{p,\mu+(j-\f r2)p}.
\end{align*}
The desired conclusion then follows from Theorem \ref{thm-6-2}.
\end{proof}

It remains to be seen if $D_{i,d+1}^r \wt f$ term in \eqref{Ball-bast-best}
can be bounded by a term that involves only $(\varphi \partial)^r f$ in
the case of $\mu > 0$.

Similar to the case of $\sph$, we can also define a  Lipschitz space
on the ball. %$\CW_p^{r,\a}(\BB^d, W_\mu)$ for $r\in \NN$ and $ \a\in [0,1)$.

\begin{defn}
For $r\in\NN$, $\a\in [0,1)$, and $1\leq p\leq \infty$,  we define
$ \CW_p^{r,\a}(\BB^d, W_\mu)$ to be the space of all functions $
f:\BB^d\to\RR$   with finite norm
\begin{align*}
\|f\|_{\mathcal{W}_p^{r,\a}(\BB^d, W_\mu)} :=\|f\|_{p,\mu}
& + \max_{1\leq i<j\leq d} \sup_{0<|\t|\leq 1} |\t|^{-\a} \|\tr_{i,j,\t}^\ell (D_{i,j}^r f )\|_{p,\mu}\\
& + \max_{1\leq i\leq d}  \sup_{0<|\t|\leq 1} |\t|^{-\a} \|\tr_{i,d+1,\t}^\ell
(D_{i,d+1}^r \wt f )\|_{L^p(\BB^d,W_{\mu-1/2})}
\end{align*}
with the usual change when $p=\infty$, where $\ell$ is a fixed positive integer, say $\ell =1$.
\end{defn}

We can also give an equivalent characterization of the space $\CW_p^{r,\a}(\BB^d, W_\mu)$
in terms of our modulus of smoothness. For the same set of parameters as in the definition
of $\CW_p^{r,\a}(\BB^d, W_\mu)$, we define a space
$$
 H_p^{r+\a}(\BB^d, W_\mu):= \left \{ f\in L^p(\BB^d, W_\mu):\quad \sup_{0<t\leq 1}
   \f{ \o_{r+\ell}(f,t)_{p,\mu}} {t^{r+\a}}<\infty \right \}.
$$

\begin{thm}\label{thm-LipBall}
Let $\mu = \frac{m-1}{2}$. If $r\in\NN$,  $1\leq p\leq \infty$, and $\a\in [0,1)$, then
$$
  \CW_p^{r,\a}(\BB^d, W_\mu) = H_p^{\a+r}(\BB^d, W_\mu) \quad \hbox{and}\quad
        \|f\|_{\CW_p^{r,\a}(\BB^d, W_\mu)} \sim \|f\|_{H_p^{\a+r}(\BB^d, W_\mu)}.
$$
\end{thm}

\begin{proof}
This follows from Theorem \ref{thm-2-5} since for $F(x,x') : = f(x)$,
$(x,x')\in \SS^{d+m-1}$ and $x \in \BB^d$, we have $\o(F,t)_{L^p(\SS^{d+m-1})}
\sim \o_{r}(f,t)_{p,\mu}$ by Lemma 5.4 of \cite{DaX}.
\end{proof}
 
We could also define a Lipschitz space that uses central differences of $\p_i^r f$ in place
of $D_{i,d+1}^r \wt f$ in $\CW_p^{r,\a}(\BB^d, W_\mu)$, so that it is equivalent to an
analogue of $H_p^{r+\a} (\BB^d, W_\mu)$ with $\o_{r+\ell}(f,t)_{p,\mu}$ in place of
$\wh \o_{r+\ell}(f,t)_{p,\mu}$.

As a consequence of the last theorem and the Jackson estimate, we have

\begin{cor} Let $\mu=\f {m-1}2$ and $m\in \NN$.
If  $r \in \NN$,  $\a\in [0,1)$,   $f \in \CW_p^{r,\a}(\BB^d,
W_\mu)$, and $1 \le p \le \infty$, then
\begin{equation} \label{LipBallestimate}
   E_n(f)_{p,\mu} \leq  c n^{-r-\a}\|f\|_{\CW_p^{r,\a}(\BB^d, W_\mu)}.
\end{equation}
\end{cor}

Let us point out that for $f \in C^r(\BB^d)$ the traditional  definition of the Lipschitz
continuity takes the form, for $0< \a < 1$,
\begin{equation}\label{Lip}
   |\p^\beta f(x) -\p^\beta f(y)| \le c \|x - y \|^\a, \quad \beta \in \NN^d, \quad |\beta| = r
\end{equation}
for all $x, y \in \BB^d$. Let us denote by $\mathrm{Lip}_{r, \a}$ the space of all $C^r(\BB^d)$
functions that satisfy \eqref{Lip}. From the definition of $D_{i,j}$ it follows readily that
$$
    \mathrm{Lip}_{r, \a} \subset \CW_\infty^{r,\a}.
$$ 
Hence, the estimate \eqref{LipBallestimate} holds for the functions in $\mathrm{Lip}_{r,\a}$. 
On the other hand, our definition of $\CW_p^{r , \a}$ is more general than $\mathrm{Lip}_{r,\a}$ 
as the following example shows.

\medskip\noindent
{\it Example.}  Let $f_\a(x) = (1-\|x\|^2 + \|x - x_0\|^2)^\a$ on $\BB^d$ with a fixed
$x_0 \in \SS^{d-1}$. Assume $1/2 < \a < 1$. Then by \cite[Example 10.1]{DaX},
$\o_r(f_\a,t)_\infty  \sim t^{2\a}$, so that by Theorem 5.8, $f_\a \in \CW_\infty^{1, 2 \a -1}(\BB^d)$.
On the other hand, setting $x_0 = (1,0,0,\ldots, 0)$ shows that $f_\a(x) =( 1-x_1^2 + (1-x_1)^2)^\a
 = 2^\a (1 - x_1)^\a$, whose first partial derivative is unbounded on $\BB^d$ so that it is not an  
$\mathrm{Lip}_{1,\a}$ function. We note that $D_{i,j} \tilde {f_\a} \in C(\BB^{d+1})$ for all 
$1\leq i<j\leq d+1$.

\end{document}